\documentclass{amsart}
\usepackage{amssymb}
\usepackage{amsmath}
\usepackage{amsthm}
\newcommand{\R}{\mathbb{R}}

\newtheorem{cor}{Corollary}
\newtheorem{theorem}{Theorem}
\newtheorem{prop}{Proposition}
\newtheorem{lemma}{Lemma}
\theoremstyle{remark}
\newtheorem*{remarks}{Remarks}
\newtheorem*{remark}{Remark}
\newtheorem*{ack}{Acknowledgements}
\usepackage{parskip}
\makeatletter
\def\thm@space@setup{%
 \thm@preskip=\parskip \thm@postskip=0pt
}
\makeatother

\title{Stability of trace theorems on the sphere}
\author{Neal Bez, Chris Jeavons, Tohru Ozawa and Mitsuru Sugimoto}
\address{Neal Bez, Department of Mathematics, Graduate School of Science and Engineering, Saitama University, Saitama 338-8570, Japan}
\email{nealbez@mail.saitama-u.ac.jp}
\address{Chris Jeavons and Tohru Ozawa, Department of Applied Physics, Waseda University, Tokyo 169-8555, Japan}
\email{iac15222@kurenai.waseda.jp, chrispjeavons1@gmail.com}
\email{txozawa@waseda.jp}
\address{Mitsuru Sugimoto, Graduate School of Mathematics, Nagoya University, Furocho, Chikusa-ku, Nagoya 464-8602, Japan}
\email{sugimoto@math.nagoya-u.ac.jp}

\begin{document}

\begin{abstract}
We prove stable versions of trace theorems on the sphere in $L^2$ with optimal constants, thus obtaining rather precise information regarding near-extremisers. We also obtain stability for the trace theorem into $L^q$ for $q > 2$, by combining a refined Hardy--Littlewood--Sobolev inequality on the sphere with a duality-stability result proved very recently by Carlen. Finally, we extend a local version of Carlen's duality theorem to establish local stability of certain Strichartz estimates for the kinetic transport equation. 
\end{abstract}
\maketitle
\section{Introduction}
For $n\geq 2$, consider the fractional Sobolev inequality for functions on $\R^{n-1}$, which for later convenience we state as
\begin{equation}\label{FS}
\|G\|_{L^{q}(\R^{n-1})}^2\leq C_{\mathrm{FS}}(n,s)\|(-\Delta)^{\frac{2s-1}{4}}G\|_{L^2(\R^{n-1})}^2,
\end{equation}
where $q=\frac{2(n-1)}{n-2s}$ and $s\in(\frac{1}{2},\frac{n}{2})$. The sharp constant and characterisation of extremisers (i.e.\ nontrivial cases of equality) for this estimate was given for $s=\frac{3}{2}$ by Aubin \cite{Aubin} and Talenti \cite{Talenti} independently, and for general $s$ by Lieb \cite{L}. Further, it is known that one may prove refinements, or stable versions, of \eqref{FS} by adding a term proportional to the distance to the set of extremisers (we denote this set by $M_{\mathrm{FS}}=M_{\mathrm{FS}}(n,s)$) to the right-hand side:
\begin{equation}\label{CFW-eq}
C_{\mathrm{FS}}(n,s)\|(-\Delta)^{\frac{2s-1}{4}}G\|_{L^2}^2-\|G\|_{L^{q}}^2\geq \alpha \inf_{G_\ast\in M_{\mathrm{FS}}}\|(-\Delta)^{\frac{2s-1}{4}}(G-G_\ast)\|_{L^2}^2,
\end{equation}
for some $\alpha>0$. Inequality \eqref{CFW-eq} was first proved for $s=\frac{3}{2}$ by Bianchi--Egnell \cite{BE}, extended to some additional values of $s$ in \cite{BWW} and \cite{LW}, and finally completed for all admissible values of $s$ recently by Chen--Frank--Weth \cite{CFW}.

The main purpose of this article is to establish stable versions of trace theorems on spheres. The classical trace theorems on the unit sphere allow us to give meaning to the restriction to $\mathbb{S}^{n-1}$ of a function defined on $\mathbb{R}^n$, assuming the function is sufficiently regular. This regularity may be captured by use of Sobolev spaces, in which case there are differences depending on whether we mean the homogeneous or inhomogeneous versions of these spaces. A natural way to unify and generalise trace theorems associated with these spaces was exposed in \cite{BSS} and we shall prove our stability theorem in the same level of generality. The underlying inequality takes the form
\begin{equation} \label{gen-trace}
\| \mathcal{S}_wg\|_{L^2(\mathbb{S}^{n-1})} \leq C(w) \|g\|_{L^2(\mathbb{R}^n)}
\end{equation}
where $\mathcal{S}_w = \mathcal{R} \circ w(\sqrt{-\Delta})^{\frac{1}{2}}$, $\mathcal{R}$ denotes the operation of restriction to $\mathbb{S}^{n-1}$, and the function $w : (0,\infty) \to (0,\infty)$ is such that the Fourier transform $\widehat{w(|\cdot|)}$ makes sense (at least away from the origin) and is positive. The boundedness of $\mathcal{R}$ from either the homogeneous Sobolev space $\dot{H}^s(\mathbb{R}^n)$, or the inhomogeneous space $H^s(\mathbb{R}^n)$, to $L^2(\mathbb{S}^{n-1})$ are the classical cases of interest. Such bounds are clearly equivalent to \eqref{gen-trace} by taking $w(r) = r^{-2s}$ and $w(r) = (1 + r^2)^{-s}$, respectively; the well-definedness and positivity of the Fourier transform are well-known in both cases.

We take $C(w)$ to mean the sharp constant in \eqref{gen-trace}, in which case we understand from \cite{BSS} that, under the above hypotheses on $w$, we have
\[
C(w)^2 = \int_0^\infty J_{\frac{n-2}{2}}(r)^2 rw(r) \, \mathrm{d}r
\]
where $J_\nu$ denotes the Bessel function of the first kind of order $\nu$. In the case where $w(r) = r^{-2s}$, the sharp form of \eqref{gen-trace} was studied earlier: the value of $C(w)$ was first found in \cite{RS}, and a characterisation of extremisers was established in \cite{BMS} and independently by Beckner in \cite{Beckner}.

The following is our first main result. In order to state it we denote by $H_k$ the space of spherical harmonics of degree $k$, that is, the space of homogeneous harmonic polynomials of degree $k$ on $\R^n$ restricted to $\mathbb{S}^{n-1}$. We take the following Fourier transform
\[
\widehat{g}(\xi)=\int_{\R^n}g(x)e^{-ix\cdot\xi}\,\mathrm{d}x,
\]
and we use $M(\mathcal{S}_w)$ to denote the set of extremisers for \eqref{gen-trace}. We define
\begin{equation}\label{lambdastar-def}
\lambda_\ast(w)=\sup_{k\geq 1}\lambda_k(w),
\end{equation}
where
\begin{equation}\label{lambda-def}
\lambda_k(w):=\int_0^\infty J_{k+\frac{n-2}{2}}(r)^2 r w(r)\,\mathrm{d}r,
\end{equation}
and we let $\mathcal{K}=\mathcal{K}(w)$ denote the (possibly empty) set of those $k\geq 1$ for which the supremum in \eqref{lambdastar-def} is attained.
\begin{theorem}\label{cor-3} 
Suppose that $w$ is as above. Then
\begin{equation}\label{ctrace-RTeq}
C(w)^2\|g\|_{L^2(\R^n)}^2 - \|\mathcal{S}_wg\|_{L^2(\mathbb{S}^{n-1})}^2 \geq  C\inf_{g_\ast\in M(\mathcal{S}_w)}\|g-g_\ast\|_{L^2(\R^n)}^2
\end{equation}
holds with 
\[
C=C'(w):=\lambda_0(w)-\lambda_\ast(w)
\]
for any $g\in L^2$, where $\lambda_k(w)$ is given by \eqref{lambda-def} and $\lambda_\ast(w)$ by \eqref{lambdastar-def}. If $w$ is such that $C'(w)>0$ then the constant is optimal, and \eqref{ctrace-RTeq} has an extremiser if and only if $\mathcal{K}\neq\emptyset$. In this case, equality holds in \eqref{ctrace-RTeq} if and only if there exists $c\in\mathbb{C}$ and $Y_k\in H_k$ such that
\[
\frac{|\xi|^{\frac{n-2}{2}}}{w(|\xi|)^\frac{1}{2}}\widehat{g}(\xi)=cJ_{\frac{n-2}{2}}(|\xi|)+\sum_{k\in\mathcal{K}}Y_k\left(\xi^\prime\right)J_{\frac{n-2}{2}+k}(|\xi|),
\]
for $\xi\in\R^n$, where $\xi^\prime:=\frac{\xi}{|\xi|}$ for $\xi\neq 0$. If $w$ is such that $C'(w)$ is equal to zero, then the estimate \eqref{ctrace-RTeq} is false for all $C>0$.
\end{theorem}
\begin{samepage}
\begin{remarks}\leavevmode\vspace{-1pt}
\begin{itemize}
\item It was proved in \cite{BSS} (see formula (1.10)) that we have the alternative representation
\begin{equation}\label{lambda-def2}
\lambda_k(w)=\frac{|\mathbb{S}^{n-2}|}{(2\pi)^n}\int_{-1}^1 F_w(1-t)P_{n,k}(t)(1-t^2)^{\frac{n-3}{2}}\,\mathrm{d}t,
\end{equation}
where $F_w$ is defined by $F_w(\frac{|\xi|^2}{2})=\widehat{w(|\cdot|)}(\xi)$. Here $P_{n,k}$ denotes the Legendre polynomial of degree $k$ in $n$ dimensions, whose definition may be found in \cite{AH} along with the fact that
\[
P_{n,k}(t)\leq 1 = P_{n,0}(t)
\]
for all $t\in [-1,1]$, $k\geq 1$ and $n\geq 2$. Since $F_w$ is positive we see that 
\[
\lambda_k(w)<\lambda_0(w)
\] 
for any $k\geq 1$, so we see that $C'(w)$ is non-negative, and in fact strictly positive if $\mathcal{K}\neq\emptyset$.
\item In the case where $\mathcal{K}=\emptyset$ we prove the constant in \eqref{ctrace-RTeq} is sharp by explicitly constructing an extremising sequence; this same argument applies when $\lambda_\ast(w)=\lambda_0(w)$, yielding the failure of \eqref{ctrace-RTeq} in this case.
\end{itemize}
\end{remarks}
\end{samepage}
We now discuss some applications of Theorem \ref{cor-3}; for convenience we define
\[
M(\mathcal{R})=w(\sqrt{-\Delta})^{\frac{1}{2}}M(\mathcal{S}_w).
\]
When $w(r)=r^{-2s}$ for $s\in(\frac{1}{2},\frac{n}{2})$, one has
\[
\lambda_k(w)=2^{1-2s}\frac{\Gamma(2s-1)\Gamma(k+\tfrac{n-2s}{2})}{\Gamma(s)^2\Gamma(k-1+\tfrac{n+2s}{2})}
\]
which implies $(\lambda_k(w))_{k\geq 1}$ is strictly decreasing (see \cite{BS}, proof of Theorem 1.6 and Lemma 5.1, also \eqref{bessel-2} below) and hence $\mathcal{K}=\{1\}$. Therefore, Theorem \ref{cor-3} immediately yields a stable version of the trace theorem for functions in the homogeneous Sobolev space $\dot{H}^s(\R^n)$, and the sharp constant may be given in closed form.
\begin{cor}\label{coro-3} 
For $s\in (\frac{1}{2},\frac{n}{2})$, the inequality
\begin{equation}\label{cltrace-RTeq}
\|\mathcal{R}\|^2\|f\|_{\dot{H}^s(\mathbb{R}^n)}^2 - \|\mathcal{R}f\|_{L^2(\mathbb{S}^{n-1})}^2 \geq  C\inf_{f_\ast\in M(\mathcal{R})}\|f-f_\ast\|_{\dot{H}^s(\mathbb{R}^n)}^2
\end{equation}
holds for any $f\in \dot{H}^s(\mathbb{R}^n)$ with constant 
\[
C=2^{1-2s}\frac{\Gamma(2s-1)}{\Gamma(s)^2}\left(\frac{\Gamma(\tfrac{n-2s}{2})}{\Gamma(\tfrac{n+2s-2}{2})}-\frac{\Gamma(\tfrac{n-2s+2}{2})}{\Gamma(\tfrac{n+2s}{2})}\right).
\]
The constant is optimal and equality holds in \eqref{cltrace-RTeq} if and only if there exists $c\in\mathbb{C}$ and $Y_1\in H_1$ such that
\[
|\xi|^{\frac{n-2}{2}+2s}\widehat{f}(\xi)=cJ_{\frac{n-2}{2}}(|\xi|)+Y_1\left(\xi'\right)J_{\frac{n}{2}}(|\xi|),
\]
for $\xi\in\R^n$.
\end{cor}
Corollary \ref{coro-3} of course implies a stable version of the homogeneous trace theorem which may be viewed as an analogue of the stable Sobolev inequality \eqref{CFW-eq}. This latter result also admits a reverse form (see e.g.\ \cite{CFW}), and a similar estimate holds for \eqref{cltrace-RTeq}; the proof of this is not difficult and will be postponed to Section \ref{FR}. 

The other case of particular interest in Theorem \ref{cor-3} arises from the choice $w(r)=(1+r^2)^{-s}$ for $s\in(\frac{1}{2},\infty)$, for which \eqref{gen-trace} is equivalent to the trace theorem for functions in the inhomogeneous Sobolev space $H^s(\R^n)$. 
\begin{cor}\label{coro-4} 
For $s\in (\frac{1}{2},\infty)$, the inequality
\begin{equation}\label{cltrace-RTeq2}
\|\mathcal{R}\|^2\|f\|_{H^s(\mathbb{R}^n)}^2 - \|\mathcal{R}f\|_{L^2(\mathbb{S}^{n-1})}^2 \geq  C\inf_{f_\ast\in M(\mathcal{R})}\|f-f_\ast\|_{H^s(\mathbb{R}^n)}^2
\end{equation}
holds for any $f\in H^s(\mathbb{R}^n)$. The constant $C$ is as in Theorem \ref{cor-3}, is positive and optimal, and there exists an extremiser for \eqref{cltrace-RTeq2}. If in addition $s=1$ then we have
\[
C=I_{\frac{n-2}{2}}(1)K_{\frac{n-2}{2}}(1)-I_{\frac{n}{2}}(1)K_{\frac{n}{2}}(1)
\]
where $I_\mu$ and $K_\mu$ are modified Bessel functions of the first kind of order $\mu$. In this case, equality holds in \eqref{cltrace-RTeq2} if and only if there exists $c\in\mathbb{C}$ and $Y_1\in H_1$ such that
\[
|\xi|^{\frac{n-2}{2}}(1+|\xi|^2)\widehat{f}(\xi)=cJ_{\frac{n-2}{2}}(|\xi|)+Y_1\left(\xi'\right)J_{\frac{n}{2}}(|\xi|),
\]
for $\xi\in\R^n$.
\end{cor}
In the case $s=1$, the proof of Corollary \ref{coro-4} is immediate from Theorem \ref{cor-3} and the observation in \cite{BS} that we may write
\begin{equation}\label{lambda-inhom}
\lambda_k(w)=I_{k+\frac{n-2}{2}}(1)K_{k+\frac{n-2}{2}}(1),
\end{equation}
which again implies $(\lambda_k(w))$ is decreasing (see \cite{BSS}). Although we do not know of such a precise formula as \eqref{lambda-inhom} for general $s$ the proof of Corollary \ref{coro-4} proceeds using a further analysis of the sequence $(\lambda_k(w))$ based on the formula \eqref{lambda-def} and some well-known estimates for the Bessel functions which arise.

Our next result concerns the homogeneous trace theorem into $L^q$. Specifically, if we fix $w(r)=r^{-2s}$ for $s\in(\frac{1}{2},\frac{n}{2})$ and $q=\frac{2(n-1)}{n-2s}>2$ then it is known that the operator $\mathcal{S}:=\mathcal{S}_w$ in fact maps into the smaller space $L^q(\mathbb{S}^{n-1})$, i.e.\ we have the inequality
\begin{equation}\label{q-trace}
\left\|\mathcal{S}g\right\|_{L^{q}(\mathbb{S}^{n-1})}\leq C_{\mathrm{Tr}}(n,s)\|g\|_{L^2(\R^n)},
\end{equation}
where the constant is taken to be optimal. The sharp inequality \eqref{q-trace} and a characterisation of extremisers\texttt{} was proved in \cite{BMS} and independently by Beckner in \cite{Beckner}; this follows from a duality argument which shows, ultimately, that \eqref{q-trace} is equivalent to the sharp Sobolev inequality \eqref{FS}. 
\begin{theorem}\label{cor-2} 
Suppose that $s\in (\frac{1}{2},\frac{n}{2})$, and $q=\frac{2(n-1)}{n-2s}$. Then\footnote{Here and throughout, we use notation $A\lesssim B$ and $A\gtrsim B$ to denote, respectively, $A\leq cB$ and $A\geq cB$ for an arbitrary constant $c>0$ which may depend on numbers such as $p$ or $q$ but never on functions such as $f$ or $G$. The value of $c$ may change from line to line.}
\begin{equation}\label{trace-RTeq}
C_{\mathrm{Tr}}(n,s)^2\|g\|_{L^2(\R^n)}^2 - \|\mathcal{S}g\|_{L^q(\mathbb{S}^{n-1})}^2 \gtrsim \inf_{f_\ast\in M(\mathcal{S})}\|g-g_\ast\|_{L^2(\R^n)}^2
\end{equation}
holds for any $g\in L^2(\mathbb{R}^n)$, where $M(\mathcal{S})$ denotes the set of extremisers for \eqref{q-trace}.
\end{theorem} 
Despite the equivalence of \eqref{FS} and \eqref{q-trace} as demonstrated in \cite{Beckner} and \cite{BMS}, it is not clear to us how to adapt the proof of \eqref{CFW-eq} from \cite{CFW} into a proof of \eqref{trace-RTeq}. One difficulty comes from the fact that the conformal invariance, which is well-known to hold for \eqref{FS} and is an important tool in the proof of \eqref{CFW-eq}, seems difficult to show for \eqref{q-trace} directly. A manifestation of this difficulty comes from the shape of the extremisers for \eqref{q-trace}; they satisfy
\[
g=|\cdot|^{s-n}\ast J_\tau^{\frac{1}{q}} \,\mathrm{d}\sigma,
\]
where $J_\tau$ is the jacobian of a conformal transformation on $\mathbb{S}^{n-1}$, and $\mathrm{d}\sigma$ is surface measure (see \cite{BMS}).  

Our proof of Theorem \ref{cor-2} instead relies on a very recent result of Carlen \cite{Carlen} which implies that the operation of `dualising' an inequality is stable under refinements such as \eqref{CFW-eq}. It is well-known that by duality, \eqref{FS} is equivalent to the following Hardy--Littlewood--Sobolev inequality, which was first established in sharp form by Lieb \cite{L}:
\begin{equation}\label{HLS}
\left|\int_{\R^{2(n-1)}}\frac{f(x)\overline{f(y)}}{|x-y|^{n-2s}}\,\mathrm{d}x\mathrm{d}y\right|\leq C_{\mathrm{HLS}}(n,s)\|f\|_{L^{q'}(\R^{n-1})}^2
\end{equation}
for $f\in L^{q'}(\R^{n-1})$, where $C_{\mathrm{HLS}}(n,s)$ is the sharp constant and $q':=\frac{q}{q-1}$ is the usual H\"older conjugate. As a consequence, in \cite{Carlen} a refinement of \eqref{HLS} is obtained using \eqref{CFW-eq}:
\begin{equation}\label{LZ-eq}
C_{\mathrm{HLS}}(n,s)\|f\|_{L^{q'}}^2-\left|\int_{\R^{2(n-1)}}\frac{f(x)\overline{f(y)}}{|x-y|^{n-2s}}\,\mathrm{d}x\mathrm{d}y\right|\gtrsim \inf_{f_\ast\in M_{\mathrm{HLS}}}\|f-f_\ast\|_{L^{q'}}^2,
\end{equation}
where $M_{\mathrm{HLS}}=M_{\mathrm{HLS}}(n,s)$ denotes the set of functions for which one has equality in \eqref{HLS}. Given \eqref{LZ-eq} and the perspective in \cite{Carlen}, a revisit of the proof of the sharp inequality \eqref{q-trace} from \cite{BMS} yields \eqref{trace-RTeq} with little additional effort.

\begin{remark}
The inequality \eqref{LZ-eq} appears to have been first proved by Liu--Zhang via a direct derivation (\cite{LZ}, Theorem 2.2), however the approach in \cite{Carlen} permits some simplifications compared to this. For instance, a lack of smoothness of the functional given by the left-hand side of \eqref{LZ-eq} occurs since $q>2$, and this causes a failure of the second-order Taylor expansion used in the derivation of \eqref{CFW-eq} (necessitating a result of Christ from \cite{C}; see Section 4.2 of \cite{LZ}).
\end{remark}

In the years since the influential work \cite{L}, there has been considerable progress made in the understanding of other inequalities with conformal structure. In order to motivate our next result we consider the kinetic transport equation
\begin{equation}\label{KTE}
\begin{cases}
& \partial_tF(t,x,v)+v\cdot\nabla_xF(t,x,v)=0 \\ 
& F(0,x,v)=f(x,v),
\end{cases}
\end{equation}
for $n\geq 1$, where $(t,x,v)\in \R\times\R^n\times\R^n$. For the velocity average
\[
\rho f(t,x):=\int_{\R^n} F(t,x,v)\,\mathrm{d}v=\int_{\R^n}f(x-tv,v)\,\mathrm{d}v,
\]
a full range of necessary and sufficient conditions on $(p,q,r)$ for the estimate
\begin{equation}\label{KT}
\|\rho f\|_{L_t^q L_x^r(\R^{n+1})}\lesssim\|f\|_{L_{x,v}^{p}(\R^{2n})}
\end{equation}
to hold is now known; see \cite{BBGL}, \cite{CP} and \cite{KT}. Since the adjoint operator is given by
\[
\rho^\ast G(x,v)=\int_\R G(s,x+vs)\,\mathrm{d}s, 
\]
then for $p=\frac{n+2}{n+1}$ and $q=r=\frac{n+2}{n}$, the dual inequality to \eqref{KT}
\begin{equation}\label{DKT}
\|\rho^\ast G\|_{L_{x,v}^{n+2}(\R^{2n})}\lesssim\|G\|_{L_{t,x}^{\frac{n+2}{2}}(\R^{n+1})}
\end{equation}
has an alternative interpretation as an estimate for the classical X-ray transform on $\R^{n+1}$; see \cite{BBI} for discussion and further references regarding these estimates (and natural generalisations such as the Radon transform) from this perspective. In particular, it follows from the results of \cite{C2}, \cite{Drouot}, and \cite{F} that the optimal constant in \eqref{DKT} is attained when
\[
G(t,x)=\frac{1}{1+t^2+|x|^2}
\]
uniquely, up to the invariances of the inequality. This result implies (cf.\ Lemma \ref{th1-i}, below) that for such $p,q$ and $r$, the optimal constant in \eqref{KT} is attained if and only if
\[
f(x,v)=\frac{1}{\left((1+|x|^2)(1+|v|^2)-(x\cdot v)^2\right)^{\frac{n+1}{2}}},
\]
again up to invariances. We remark that the proof of the sharp form of \eqref{DKT} relies on the conformal invariance of the inequality in this case; for general triples $(p,q,r)$ this fails and it is an open problem to obtain the sharp constant and characterisation of extremisers in \eqref{KT}.

Developing the ideas of \cite{C2}, \cite{Drouot}, and \cite{F} further, in \cite{Drouot2} it is proved that the inequality \eqref{DKT} is \emph{locally stable} in the following sense: the inequality
\begin{equation}\label{D-eq}
\|\rho\| - \frac{\|\rho^\ast G\|_{L^{p'}}}{\|G\|_{L^{q'}}}\gtrsim \inf_{G_\ast\in M(\rho^\ast)}\left(\frac{\|G-G_\ast\|_{L^{q'}}}{\|G\|_{L^{q'}}}\right)^2
\end{equation}
holds for any $G\in L^{q'}\setminus\{0\}$ such that
\[
\inf_{G_\ast\in M(\rho^\ast)}\frac{\|G-G_\ast\|_{L^{q'}}}{\|G\|_{L^{q'}}}\leq c_1
\] 
for some computable constant $c_1 < 1$ depending only on $n$. Here, as usual $M(\rho^\ast)$ denotes the set of extremisers for \eqref{DKT}. The notion of local stability has proved useful more generally: such estimates are more tractable than the corresponding global ones but turn out to be equivalent in many cases (for example, this holds for \eqref{FS} and is used in the proof of \eqref{CFW-eq}; see \cite{CFW}). Our next result is the following, which implies that the inequality \eqref{KT} is locally stable.
\begin{theorem}\label{cor-4}
For $p=\frac{n+2}{n+1}$ and $q=r=\frac{n+2}{n}$,
\begin{equation}\label{KT-RTeq}
\|\rho\| - \frac{\|\rho f\|_{L^q}}{\|f\|_{L^p}} \gtrsim \inf_{f_\ast\in M(\rho)}\left(\frac{\|f-f_\ast\|_{L^{p}}}{\|f\|_{L^{p}}}\right)^2
\end{equation}
for $f\in L^p(\R^{2n})\setminus\{0\}$ such that
\[
\inf_{f_\ast\in M(\rho)}\frac{\|f-f_\ast\|_{L^{p}}}{\|f\|_{L^{p}}}\leq c_2,
\]
where $c_2<1$ is a computable constant. Here $M(\rho)$ denotes the set of extremisers for \eqref{KT} for such $(p,q,r)$.
\end{theorem}
\begin{samepage}
\begin{remarks}\leavevmode\vspace{-1pt}
\begin{itemize}
\item Although it is natural to expect \eqref{D-eq} to hold globally this is not known for $n\geq 2$. As noted in \cite{Drouot2}, it would be enough to establish that any extremising sequence for \eqref{DKT} is precompact, but the inequality has a particularly large group of symmetries which leads to many ways in which this can fail. It is known when $n=1$ (due to Christ - see \cite{C2}) and so the global estimate holds in this case, but for $n=1$ the inequality \eqref{KT} is self-dual and so \eqref{D-eq} and \eqref{KT-RTeq} are the same.
\item Corollary \ref{cor-4} is related to a number of recent results where properties of the solutions to \eqref{KTE} have been studied using the relation with the X-ray transform outlined above. Examples include monotonicity properties for estimates closely related to \eqref{DKT} (\cite{BBI}, Section 5), and null-form estimates which may be viewed as multilinear generalisations of \eqref{KT} (\cite{BBGL2}, Section 7.2).
\end{itemize}
\end{remarks}
\end{samepage}
Our proof of Theorem \ref{cor-4} proceeds using \eqref{D-eq} and a duality theorem for local stability inequalities which is of independent interest; see Theorem \ref{t1-loc}, below. It is inspired by, and is a generalisation of, a result from \cite{Carlen} which was used to prove a local version of \eqref{LZ-eq} in the case $s=\frac{3}{2}$. 

The question of stability has been studied recently for a number of important geometric and analytic inequalities and has found a range of applications (see \cite{CF}, \cite{CFL}, \cite{CS}). However, the derivation of estimates of this type with optimal distance norms has in general proved to be a difficult problem (see e.g.\ \cite{BJ}, \cite{C-Young}, \cite{C}, \cite{CFMP}, \cite{FMP} and references therein) and many important questions remain open. Our results are elementary in comparison but as far as we know, Theorem \ref{cor-3} provides the first example of a global Bianchi--Egnell type stability estimate with optimal constant. Further, Theorems \ref{cor-2} and \ref{cor-4} build on the very recent work \cite{Carlen}, using this perspective to derive estimates for which we are not aware of an alternative approach.

\emph{Organisation.} The next section is devoted to the results relating to stability of trace theorems into $L^2$; specifically we prove Theorem \ref{cor-3} and show how to deduce Corollary \ref{coro-4}. In Section \ref{s-cor} we prove the stable $L^q$ trace inequality Theorem \ref{cor-2}, and in Section \ref{s-cor4} we state and prove the local duality-stability result that we need to deduce Theorem \ref{cor-4}. Finally in Section \ref{FR} we prove the reverse form of \eqref{cltrace-RTeq} mentioned above and discuss some related problems.

\section{Stability of trace estimates into $L^2$ - proof of Theorem \ref{cor-3} and Corollary \ref{coro-4}}
Our proof of Theorem \ref{cor-3} is based on a decomposition of $L^2(\R^n)$ induced by the spherical harmonics. In order to understand sharp smoothing estimates for linear dispersive propagators, Walther \cite{Walther} used such a decomposition (later developed in \cite{BS} and \cite{BS2}); these smoothing estimates are known to imply trace theorems on the sphere and thus our approach is natural. In particular we use the following decomposition of $L^2(\R^n)$ (\cite{BS}, Section 2.1):
\begin{equation}\label{f-dec}
\widehat{g}(\xi)=\sum_{k\geq 0}\sum_{m=1}^{\mathrm{dim}(H_k)}P^{(k,m)}\left(\frac{\xi}{|\xi|}\right)g_{0}^{(k,m)}(|\xi|)|\xi|^{\frac{1-n}{2}}, \;\; \xi\in\R^n,
\end{equation}
where for each $k$ the collection $\{P^{(k,m)}\}_{m=1}^{\mathrm{dim}(H_k)}$ forms an orthonormal basis of the space of spherical harmonics $H_k$, and $g_{0}^{(k,m)}\in L^2(0,\infty)$ for each $k\geq 0$ and $1\leq m\leq \mathrm{dim}(H_k)$. Unless otherwise specified, by $(k,m)$ we mean pairs of integers $k$ and $m$ in this range, and we use notation $\sum_{k,m}$ for the double sum in \eqref{f-dec}. For any $P_k\in H_k$ one has (\cite{Walther}, Corollary 5.1)
\[
\widehat{P_k\mathrm{d}\sigma}(x)=\frac{(2\pi)^{\frac{n}{2}}}{i^k}P_k\left(\frac{x}{|x|}\right)|x|^{\frac{2-n}{2}}J_{k+\frac{n-2}{2}}(|x|)
\]
for $x\in\R^n$, where $\mathrm{d}\sigma$ denotes induced Lebesgue measure on $\mathbb{S}^{n-1}$. By Fourier inversion, \eqref{f-dec} and polar co-ordinates we deduce that
\[
w(\sqrt{-\Delta})^{\frac{1}{2}}g(\theta)=\frac{1}{(2\pi)^{\frac{n}{2}}}\sum_{k,m}(-1)^k\frac{P^{(k,m)}(\theta)}{i^k}\int_0^\infty g_0^{(k,m)}(r)r^{\frac{1}{2}}w(r)^{\frac{1}{2}}J_{k+\frac{n-2}{2}}(r)\,\mathrm{d}r,
\]
for $\theta\in\mathbb{S}^{n-1}$. Hence, by orthogonality,
\begin{align*}
(2\pi)^n\|\mathcal{S}_wg\|_{L^2(\mathbb{S}^{n-1})}^2=\sum_{k,m}\left(\int_0^\infty g_0^{(k,m)}(r)r^{\frac{1}{2}}J_{k+\frac{n-2}{2}}(r)w(r)^{\frac{1}{2}}\,\mathrm{d}r\right)^2=:\sum_{k,m}A_{k,m}.
\end{align*}
Define
\[
B_{k,m}=\int_0^\infty |g_0^{(k,m)}(r)|^2\,\mathrm{d}r.
\]
Recalling that
\[
\lambda_k(w)=\int_0^\infty J_{k+\frac{n-2}{2}}(r)^2rw(r)\,\mathrm{d}r,
\]
by the Cauchy--Schwarz inequality we have $A_{k,m}(w)\leq \lambda_k(w)B_{k,m}$ for each $(k,m)$ with equality if and only if there exists constants $c_{k,m}\in\mathbb{C}$ such that
\begin{equation}\label{eq-c-1}
g_0^{(k,m)}(r)=c_{k,m}J_{k+\frac{n-2}{2}}(r)w(r)^{\frac{1}{2}}r^{\frac{1}{2}}
\end{equation}
almost everywhere on $(0,\infty)$. Further, using Plancherel's theorem and orthogonality, one has
\[
\|g\|_{L^2(\R^n)}^2=\frac{1}{(2\pi)^n}\sum_{k,m}B_{k,m},
\]
from which it follows that $C(w)^2=\lambda_0(w)$ and the image under the Fourier transform of $M(\mathcal{S}_w)$ is given by
\[
\widehat{M(\mathcal{S}_w)}=\left\{\mu w(|\cdot|)^{\frac{1}{2}}\frac{J_{\frac{n-2}{2}}(|\cdot|)}{|\cdot|^{\frac{n-2}{2}}}:\mu\in\mathbb{C}\setminus\{0\}\right\}
\]
(see Theorem 1.8 of \cite{BSS} for the original derivation of this). Using this and the Hilbert structure of $L^2$ we may assume that
\[
\inf_{g_\ast\in M(\mathcal{S}_w)}\|g-g_\ast\|_{L^2}^2=\|g\|_{L^2}^2-\frac{|\left\langle g,g_\ast\right\rangle|^2}{\|g_\ast\|_{L^2}^2},
\]
where (by abuse of notation) $g_{\ast}\in M(\mathcal{S}_w)$ is fixed. Using polar co-ordinates, orthogonality and the fact that $P^{(0,1)}$ is constant, one can compute that
\[
\frac{|\left\langle g,g_\ast\right\rangle|^2}{\|g_\ast\|_{L^2}^2}=\frac{|\mathbb{S}^{n-1}|}{(2\pi)^n\lambda_0(w)}(P^{(0,1)})^2A_{0,1}(w)=\frac{A_{0,1}(w)}{(2\pi)^n\lambda_0(w)}, 
\]   
where the second equality follows from the condition $\|P^{(0,1)}\|_{L^2(\mathbb{S}^{n-1})}^2=1$. Combining all of the above, we see that \eqref{ctrace-RTeq} (with $C=\lambda_0(w)-\lambda_\ast(w)$) is equivalent to
\begin{equation}\label{rd-eq}
\lambda_0(w)\sum_{k,m}B_{k,m}-\sum_{k,m}A_{k,m}(w)\geq (\lambda_0(w)-\lambda_\ast(w))\bigg(\bigg(\sum_{k,m}B_{k,m}\bigg)-\frac{A_{0,1}(w)}{\lambda_0(w)}\bigg),
\end{equation}
or
\[
\lambda_\ast(w)\sum_{k,m}B_{k,m}\geq A_{0,1}(w)\frac{\lambda_\ast(w)}{\lambda_0(w)}+\sum_{\substack{k,m\\k\geq 1}}A_{k,m}(w).
\]
But
\begin{align}\label{eq-ch} 
A_{0,1}(w)\frac{\lambda_\ast(w)}{\lambda_0(w)}+\sum_{\substack{k,m\\k\geq 1}}A_{k,m}(w)\leq B_{0,1}\lambda_\ast(w)+\sum_{\substack{k,m\\k\geq 1}}\lambda_k(w) B_{k,m} \leq \lambda_\ast(w) \sum_{k,m}B_{k,m},
\end{align}
and so \eqref{rd-eq} holds.

Equality holds in the first inequality in \eqref{eq-ch} if and only if \eqref{eq-c-1} holds for each $(k,m)$, and since $\lambda_k(w)>0$ equality holds in the second inequality only if $B_{k,m}=0$ for all $(k,m)$ with $k\notin \mathcal{K}$. Combining these two conditions and using the fact that the collection $\{P^{(k,m)}\}$ is a basis for the space $H_k$, we obtain the claimed equality condition for \eqref{ctrace-RTeq} in the case $\mathcal{K}\neq\emptyset$.

It remains to demonstrate the sharpness of the constant when $\mathcal{K}=\emptyset$. In this case, by the definition of $\lambda_\ast(w)$ we can find a subsequence $(\lambda_{h(l)}(w))$ of $(\lambda_l(w))$ with $\lambda_{h(l)}(w)\to\lambda_\ast(w)$ as $l\to\infty$. But then if we define $(g_l)\subset L^2(\R^n)$ by
\[
\frac{|\xi|^{\frac{n-2}{2}}}{(w(|\xi|))^\frac{1}{2}}\widehat{g_l}(\xi)=P^{(h(l),1)}\left(\xi^\prime\right)J_{\frac{n-2}{2}+h(l)}(|\xi|),
\]
then $A_{k,m}(w)=\lambda_k(w)B_{k,m}$ for all $(k,m)$, with both sides equal to zero if and only if $(k,m)\neq(h(l),1)$. But then testing \eqref{ctrace-RTeq} on such $g_l$ we see
\begin{align*}
\frac{\lambda_0(w)\sum_{k,m}B_{k,m}-\sum_{k,m}A_{k,m}(w)}{(\sum_{k,m}B_{k,m})-\frac{A_{0,1}(w)}{\lambda_0(w)}} & =\frac{\sum_{k,m}B_{k,m}(\lambda_0(w)-\lambda_k(w))}{\sum_{k,m}B_{k,m}}\\
& =\lambda_0(w)-\lambda_{h(l)}(w)\\
& \to \lambda_0(w)-\lambda_\ast(w),
\end{align*}
where the second equality follows from the fact that the sums have only one term, since $\lambda_k(w)>0$. This completes the proof of Theorem \ref{cor-3}.

\begin{remarks}\leavevmode\vspace{-1pt}
\begin{itemize}
\item We could have combined the results of \cite{Carlen} with the analysis of the dual inequality to \eqref{gen-trace} in \cite{BMS} and \cite{BSS} to obtain a simpler proof of \eqref{ctrace-RTeq}, but doing so does not seem to yield the additional information given in Theorem \ref{cor-3}. 
\item As far as we know, our argument gives, in particular, the first direct proof of the sharp trace theorem on the sphere \eqref{gen-trace}; all previous proofs of this result relied upon duality. We also point out that simple modifications of our argument yield similar direct proofs of the trace theorem with angular regularity (\cite{BMS}, Corollary 3.3), as well as a stable version of this estimate analogous to \eqref{ctrace-RTeq}.
\end{itemize}
\end{remarks}
We now turn to the proof of Corollary \ref{coro-4}; for the rest of this section we fix $w(r)=(1+r^2)^{-s}$ and $s\in(\frac{1}{2},\infty)$.  Before proceeding we record the following two facts about the Bessel function $J_{\nu(k)}$ for $\nu(k):=k+\frac{n-2}{2}$:
\begin{equation}\label{bessel-1}
|J_{\nu(k)}(r)|\lesssim r^{-\frac{1}{3}}
\end{equation}
for any $r>0$ (with implicit constant independent of $k$ and $n$), and
\begin{equation}\label{bessel-2}
\int_{0}^\infty |J_{\nu(k)}(r)|^2r^{1-\tau}\,\mathrm{d}r=2^{1-\tau}\frac{\Gamma(\tau-1)\Gamma(k+\frac{n-\tau}{2})}{\Gamma(\frac{\tau}{2})^2\Gamma(k+\frac{n+\tau}{2})}
\end{equation}
for any $\tau>1$. Inequality \eqref{bessel-1} is due to Landau \cite{Landau}, while the identity \eqref{bessel-2} goes back to Watson \cite{Watson}. By Stirling's formula, the left hand side of \eqref{bessel-2} converges to zero as $k\rightarrow\infty$, for any fixed $\tau>1$.

\proof[Proof of Corollary \ref{coro-4}] We will show that for this choice of $w$ the sequence $(\lambda_k(w))$ converges to zero; this suffices as it then follows that $\mathcal{K}\neq \emptyset$ (for if not then we could find a subsequence converging to $\lambda_\ast(w)>0$) so by the first remark after Theorem \ref{cor-3} we can conclude for $s\neq 1$, and the case $s=1$ follows from \eqref{lambda-inhom} as discussed in the introduction. 

We introduce two numbers $p>1$ and $\varepsilon>0$ depending on $s$, whose values will be specified below. Using \eqref{bessel-1} and H\"older's inequality we have
\begin{align*}
\int_0^\infty |J_{\nu(k)}(r)|^2rw(r)\,\mathrm{d}r & = \int_0^\infty |J_{\nu(k)}(r)|^{\frac{2}{p}+\frac{2}{p'}}rw(r)\,\mathrm{d}r \\
& \lesssim \int_0^\infty |J_{\nu(k)}(r)|^{\frac{2}{p}}r^{1-\frac{2}{3p'}}w(r)\,\mathrm{d}r, \\
& \leq \left(\int_0^\infty |J_{\nu(k)}(r)|^{2}r^{-\varepsilon}\,\mathrm{d}r\right)^{\frac{1}{p}}\left(\int_0^\infty w(r)^{p'}r^{p'-\frac{2}{3}+\varepsilon(p'-1)}\,\mathrm{d}r\right)^{\frac{1}{p'}}.
\end{align*}
Since by \eqref{bessel-2} the first term converges to zero as $k\rightarrow\infty$ for any fixed $\varepsilon>0$ (and since $\lambda_k(w)>0$ for all $k$), it is enough to find $\varepsilon>0$ and $p>1$ such that second term is finite, for our choice of $w$. To see that this is possible, by a simple change of variables we require
\[
\int_0^\infty (1+u)^{-p's}u^{\frac{1}{2}(p'-\frac{5}{3}+\varepsilon(p'-1))}\,\mathrm{d}u<\infty,
\]
but elementary considerations show that this holds whenever
\[
s>\frac{1}{2}+\frac{1}{6p'}+\frac{\varepsilon}{2p}.
\]
Recalling that $s>\frac{1}{2}$, the conclusion follows by choosing $p$ sufficiently close to 1 and $\varepsilon$ sufficiently small.\endproof

\begin{remark}
For $s\in\{\frac{n-1}{2},\frac{n+1}{2}\}$, it is known that there is a constant $c>0$ (depending on $n$ and $s$) such that
\begin{equation}\label{spec-cases}
\widehat{w(|\cdot|)}=c|\cdot|^{s-(\frac{n+1}{2})}e^{-|\cdot|}.
\end{equation}
When $s=\frac{n+1}{2}$, \eqref{spec-cases} is an elementary fact about the Poisson kernel. When $s=\frac{n-1}{2}$, it is likely to be known although we could not find it explicitly in the literature; it may be verified directly for $n\in\{2,3\}$ (\cite{Foschi}, Lemma 5.3 and Lemma 6.3, taking there $t=0$), from which the general case follows by an induction argument (\cite{BJOS}, proof of Theorem 2.2). It is conceivable that a more explicit value of the constant $C$ in \eqref{cltrace-RTeq2} may be obtained in these cases using \eqref{spec-cases} and the alternative formula \eqref{lambda-def2} for $\lambda_k(w)$, but we do not attempt this here.
\end{remark}

\section{Stability of trace estimates into $L^q$ - proof of Theorem \ref{cor-2}}\label{s-cor}
Our proof of Theorem \ref{cor-2} rests on the following result, which is contained in Theorem 3.3 of \cite{Carlen}. In order to state it we define, for a bounded linear operator $T$ from $L^p(X)$ to $L^q(Y)$, the set of extremisers
\[
M(T)=\{g\in L^p(X)\setminus\{0\}: \|Tg\|_{L^q(Y)}=\|T\|\|g\|_{L^p(X)}\}.
\]
In what follows, we will shorten $\|\cdot\|_{L^r}=\|\cdot\|_r$ where there is no chance of confusion.
\begin{theorem}\label{main-thm}
Let $T$ be as above for $1<p\leq 2$ and $1<q<\infty$, and such that $M(T)\neq\emptyset$. If
\begin{equation}\label{Tstar-rteq}
\|T\|^2\|G\|_{q'}^2 - \|T^\ast G\|_{p'}^2 \gtrsim \inf_{G_\ast\in M(T^\ast)}\|G-G_\ast\|_{q'}^{2}
\end{equation}
holds for $G\in L^{q'}$, then
\begin{equation}\label{T-rteq}
\|T\|^2\|g\|_p^2 - \|Tg\|_{q}^2 \gtrsim \inf_{g_\ast\in M(T)}\|g-g_\ast\|_{p}^{2}
\end{equation}
for any $g\in L^p$.
\end{theorem}
For the rest of this section, we fix $s\in(\frac{1}{2},\frac{n}{2})$, $p=2$ and $q=\frac{2(n-1)}{n-2s}$, and we recall the definition $\mathcal{S}=\mathcal{R}(-\Delta)^{-\frac{s}{2}}$. In view of Theorem \ref{main-thm}, it is enough to prove \eqref{Tstar-rteq} for $T=\mathcal{S}$. In order to do this we use notation from \cite{CFW}, introducing the isometric map $\mathcal{P}:L^{q'}(\mathbb{S}^{n-1})\rightarrow L^{q'}(\R^{n-1})$ defined by 
\[
\mathcal{P}G(x):=J_{\pi^{-1}}(x)^{\frac{1}{q'}}G(\pi^{-1}x), \qquad x\in\R^{n-1},
\]
where $\pi^{-1}:\mathbb{R}^{n-1}\rightarrow \mathbb{S}^{n-1}$ is the inverse stereographic projection
\[
\pi^{-1}(x):=\left(\frac{2x}{1+|x|^2}, \frac{1-|x|^2}{1+|x|^2}\right), \qquad x\in\R^{n-1}
\]
and $J_{\pi^{-1}}$ is the jacobian of this transformation, which one can check is
\[
J_{\pi^{-1}}(x)=\left(\frac{2}{1+|x|^2}\right)^{n-1},
\] 
again for $x\in\R^{n-1}$. We then have that
\begin{align*}
\|\mathcal{S}^{\ast}G\|_{L^2(\R^n)}^2 & =C_{n,s}\left|\int_{\mathbb{S}^{n-1}}\int_{\mathbb{S}^{n-1}}G(\omega)\overline{G(\eta)}(1-\omega\cdot\eta)^{s-\frac{n}{2}}\,\mathrm{d}\sigma_\omega \mathrm{d}\sigma_\eta\right|\\
& = C_{n,s}\left|\int_{\R^{n-1}}\int_{\R^{n-1}}\frac{\mathcal{P}G(x)\overline{\mathcal{P}G(y)}}{|x-y|^{n-2s}}\,\mathrm{d}x\mathrm{d}y\right|,
\end{align*}
where $C_{n,s}$ is a constant. The first inequality here is proved in \cite{BMS}, while the second is a well-known fact about the inequality \eqref{HLS} and goes back to \cite{L}. The constant $C_{n,s}$ is explicitly computable but we do not record its value here, referring instead to \cite{BMS}. It follows that
\begin{equation}\label{S2HLS}
C_{\mathrm{Tr}}(n,s)^2= C_{n,s}C_{\mathrm{HLS}}(n,s), \;\; \mathcal{P}M(\mathcal{S}^\ast)=M_{\mathrm{HLS}}(n,s),
\end{equation}
which in view of Lieb's result \cite{L} leads to the sharp inequality
\[
\|\mathcal{S}^\ast G\|_{L^2(\mathbb{R}^n)}\leq C_{\mathrm{Tr}}(n,s)\|G\|_{L^{q'}(\mathbb{S}^{n-1})}
\]
and characterisation of $M(\mathcal{S}^\ast)$ as proved in \cite{BMS}. Applying \eqref{LZ-eq} and using \eqref{S2HLS}, it follows that
\begin{align*}
C_{\mathrm{Tr}}(n,s)^2\|G\|_{L^{q'}(\mathbb{S}^{n-1})}^2-\|\mathcal{S}^{\ast}G\|_{L^2(\R^n)}^2 & \gtrsim \inf_{G_\ast\in M(\mathcal{S}^\ast)}\|\mathcal{P}G-\mathcal{P}G_\ast\|_{L^{q'}(\R^{n-1})}^2 \\
& = \inf_{G_\ast\in M(\mathcal{S}^\ast)}\|G-G_\ast\|_{L^{q'}(\mathbb{S}^{n-1})}^2,
\end{align*}
which is \eqref{Tstar-rteq} for $T=\mathcal{S}$.

\begin{remark}  
As noted in \cite{Carlen}, it is possible to bound the implicit constant in \eqref{T-rteq} from below in terms of the implicit constant in \eqref{Tstar-rteq}. As such, the above argument shows that the implicit constant in \eqref{trace-RTeq} is controlled by the one in \eqref{LZ-eq}, which in turn admits a bound in terms of $\alpha$ in \eqref{CFW-eq} (see \cite{Carlen}, Theorem 1.1).
\end{remark}

\section{Local duality-stability - proof of Theorem \ref{cor-4}}\label{s-cor4}
Theorem \ref{cor-4} will follow immediately from the following result which may be viewed as an analogue of Theorem \ref{main-thm} for local stability estimates. It is closely related to the result in Section 3.2 of \cite{Carlen} in which a local version of \eqref{LZ-eq} in the case $s=\frac{3}{2}$ is obtained using an appropriate local version of \eqref{CFW-eq}; we generalise this argument to a framework suitable for our application by using the results of \cite{CFL}.
\begin{theorem}\label{t1-loc}
For $1<p\leq 2$ and $1<q<\infty$, let $T:L^p(X)\rightarrow L^q(Y)$ be a bounded linear operator. Assume that 
\begin{equation}\label{Tstar-rteq2}
\|T\| - \frac{\|T^\ast G\|_{p'}}{\|G\|_{q'}} \gtrsim \inf_{G_\ast\in M(T^\ast)}\left(\frac{\|G-G_\ast\|_{q'}}{\|G\|_{q'}}\right)^2
\end{equation}
holds for any $G\in L^{q'}\setminus\{0\}$ with
\[
\inf_{G_\ast\in M(T^\ast)}\frac{\|G-G_\ast\|_{q'}}{\|G\|_{q'}}< c_1
\] 
for some $c_1\leq 1$. There is a constant $c_2\leq 1$ (with explicit dependence on $p,q,\|T\|$ and $c_1$) such that if $g\in L^p\setminus\{0\}$ is such that
\begin{equation}\label{f-near}
\inf_{f_\ast\in M(T)}\frac{\|g-g_\ast\|_p}{\|g\|_p}<c_2,
\end{equation}
then
\begin{equation}\label{T-rteq2}
\|T\| - \frac{\|Tg\|_{q}}{\|g\|_{p}} \gtrsim \inf_{g_\ast\in M(T)}\left(\frac{\|g-g_\ast\|_p}{\|g\|_p}\right)^2
\end{equation}
holds.
\end{theorem}
An important tool in the proof of this result is the following lemma, which was made explicit as a special case of Theorem 2.3 (see also Example 2.1) of \cite{Carlen} and is implicit for a number of specific operators in earlier works (see e.g.\ \cite{BMS}, \cite{CL}, \cite{C}, \cite{L}).
\begin{lemma}\label{th1-i}
Suppose that $T:L^p(X)\rightarrow L^q(Y)$ is a bounded linear operator and $1<p,q<\infty$. Then $M(T)\neq \emptyset$ if and only if $M(T^\ast)\neq\emptyset$. In this case,
\[
M(T)=|T^\ast M(T^\ast)|^{p^\prime-2}T^\ast M(T^\ast).
\]
\end{lemma}
We also use the following two inequalities from \cite{CFL}: for $r\geq 1$, if $g_1,g_2\in L^{r}$ then
\begin{equation}\label{CFL-3}
\|D_{r} g_1-D_{r} g_2\|_{r'}\leq C_r\left(\frac{\|g_1-g_2\|_{r}}{\|g_1\|_{r}+\|g_2\|_{r}}\right)^{\min\{r,2\}-1}
\end{equation}
holds with constant
\[
C_r=\begin{cases} 2(r')^{r-1} & \text{ if } r\leq 2 \\ 4(r-1) & \text{ if } r\geq 2 \end{cases}
\]
and if in addition $r\geq 2$ and $h_1$ and $h_2$ are unit vectors in $L^r$ and $L^{r'}$ respectively, then
\begin{equation}\label{CFL-1}
\left|\int h_1 h_2 \right|\leq 1-\frac{r'-1}{4}\|D_rh_1-e^{i\theta}h_2\|_{r'}^{2},
\end{equation}
where $\theta$ is such that $e^{i\theta}\int h_1h_2\geq 0$. Here, $D_r$ is the duality map
\[
D_r F:=\frac{|F|^{r-2}\overline{F}}{\|F\|_{L^r}^{r-1}},
\]
for $0\neq F\in L^r$, and we note that by homogeneity $D_r(\lambda F)=D_rF$ for all $\lambda>0$.

\proof[Proof of Theorem \ref{t1-loc}] By homogeneity it is enough to prove \eqref{T-rteq2} for $\|g\|_p=1$. Set $G=\overline{D_q Tg}$, then
\begin{align*}
\|T\|-\|Tg\|_{q} & =\|T\|-\left|\int Tg \overline{G}\right| \\
& = \|T\|-\|T^\ast G\|_{p'}\left|\int g \frac{\overline{T^\ast G}}{\|T^\ast G\|_{p'}}\right|.
\end{align*}
Taking $(h_1,h_2,r)=\left(\frac{\overline{T^\ast G}}{\|T^\ast G\|_{p'}}, g,p'\right)$ in \eqref{CFL-1} and noting that in this case $\theta=0$, we get that
\begin{align}\label{eq1}
\|T\|- \|Tg\|_{q} & \geq \|T\|-\|T^\ast G\|_{p'}\left(1-\frac{p-1}{4}\left\|g-D_{p'}\overline{T^\ast G}\right\|_{p}^2\right)\nonumber\\
& = \|T\|-\|T^\ast G\|_{p'}+\frac{p-1}{4}\|T^\ast G\|_{p'}\left\|g-D_{p'}\overline{T^\ast G}\right\|_{p}^2 \\
& \gtrsim \inf_{G_\ast}\|G-G_\ast\|_{q'}^{2}+\|T^\ast G\|_{p'}\left\|g-D_{p'}\overline{T^\ast G}\right\|_{p}^2, \nonumber
\end{align}
where we have used the homogeneity of $D_{p'}$, the hypothesis \eqref{Tstar-rteq2} and the fact that $\|G\|_{q'}=1$. Taking $(g_1,g_2,r)=(\overline{T^\ast G},\overline{T^\ast G_\ast},p')$ in \eqref{CFL-3} it follows that
\begin{align*}
\inf_{G_\ast\in M(T^\ast)}\|D_{p'}\overline{T^\ast G}-D_{p'}\overline{T^\ast G_\ast}\|_{p}^2 & \lesssim\inf_{G_\ast}\left(\frac{\|T^\ast G-T^\ast G_\ast\|_{p'}}{\|T^\ast G\|_{p'}+\|T^\ast G_\ast\|_{p'}}\right)^{2}\\
& \lesssim \left(\frac{1}{\|T^\ast G\|_{p'}}\inf_{G_\ast}\|G-G_\ast\|_{q'}\right)^{2}.
\end{align*}
Inserting this into \eqref{eq1} we see that 
\[
\|T\|- \|Tg\|_{q} \gtrsim \|T^\ast G\|_{p'}\|g-D_{p'}\overline{T^\ast G}\|_p^2+\inf_{g_\ast\in M(T)}\|T^\ast G\|_{p'}^{2}\|D_{p'}\overline{T^\ast G}-g_\ast\|_p^2,
\]
since, by Lemma \ref{th1-i}, $D_{p'}\overline{T^\ast G_\ast}\in M(T)$ whenever $G_\ast\in M(T^\ast)$. Next note that if, for example, \eqref{f-near} holds with $c_2=\frac{1}{4}$ then we can find $g_\ast\in M(T)$ with $\|g-g_\ast\|_{p}<\frac{\|g\|_p}{4}$. But then, 
\begin{equation}\label{ec1}
\|T\|-\|Tg\|_{q} =\big|\|T\|- \|T\|\|g_\ast\|_p+\|Tg_\ast\|_{q}-\|Tg\|_{q}\big| \leq 2\|T\|\|g-g_\ast\|_{p}\leq \frac{\|T\|}{2},
\end{equation}
and hence by H\"older's inequality we have that $\|T^\ast G\|_{p'}\gtrsim 1$. 

It now remains to note that if $g\in L^p$ with $\|g\|_p=1$ and $g_\ast\in M(T)$ are arbitrary, 
\begin{align}\label{ec2}
\|g-g_\ast\|_p\gtrsim \|Tg-Tg_\ast\|_{q}\gtrsim \|G-G_\ast\|_{q'}^{\eta}(\|Tg\|_q+\|Tg_\ast\|_q)\gtrsim \|G-G_\ast\|_{q'}^\eta,
\end{align}
where $G$ is as above, $G_\ast:=\overline{D_q Tg_\ast}\in M(T^\ast)$, and $\eta:=(\min\{q,2\}-1)^{-1}$. The second inequality here follows from \eqref{CFL-3}, and the third follows from \eqref{ec1}. By raising \eqref{ec2} to an appropriate power and keeping track of the constants, one may obtain an explicit dependence of $c_2$ on $c_1$ and the other quantities as claimed.\endproof

\emph{Remark.} By following the compactness argument from \cite{BE}, it is possible to use Theorem \ref{t1-loc} to give a proof of Theorem \ref{main-thm} under the assumption that
\[
\lim_{m\rightarrow\infty}\frac{\|Tg_m\|_{q}}{\|g_m\|_{p}}= \|T\| \Rightarrow \lim_{m\rightarrow\infty}\inf_{g_\ast\in M(T)}\frac{\|g_m-g_\ast\|_p}{\|g_m\|_p}= 0
\]  
for any $(g_m)\subset L^p\setminus\{0\}$. Although it turns out that this condition holds fairly generically (in particular, in view of Lemma 1 of \cite{CFW} it is enough to deduce \eqref{LZ-eq} from \eqref{CFW-eq}), proceeding in this fashion one loses information about the explicit constant in \eqref{T-rteq}, and as evidenced by Theorem 3.2 of \cite{Carlen} this need not be the case.

\section{Further results}\label{FR}

First, we establish the reverse form of \eqref{cltrace-RTeq} mentioned in the introduction. A similar result may be proved for the general estimate \eqref{ctrace-RTeq} by the same argument, using the relevant results in \cite{BSS}.
\begin{prop}
For $s\in(\frac{1}{2},\frac{n}{2})$ the inequality
\[
\|\mathcal{R}\|^2\|f\|_{\dot{H}^s(\R^n)}^2 - \|\mathcal{R}f\|_{L^2(\mathbb{S}^{n-1})}^2\leq \|\mathcal{R}\|^2\inf_{f_\ast\in M(\mathcal{R})}\|f-f_\ast\|_{\dot{H}^s(\R^n)}^2
\]
holds for any $f\in\dot{H}^s(\R^n)$.
\end{prop}
\proof We proceed following the proof from \cite{CFW} of the forthcoming inequality \eqref{CFW-eq2}. Let $f_\ast$ denote the closest point in $M(\mathcal{R})$ to $f$, then $v:=f-f_\ast$ satisfies $\langle v,f_\ast\rangle_{\dot{H}^s}=0$. We then have
\begin{align*}
\|\mathcal{R}\|^2\|f\|_{\dot{H}^s}^2 - \|\mathcal{R}f\|_{L^2(\mathbb{S}^{n-1})}^2 & =\|\mathcal{R}\|^2(\|v\|_{\dot{H}^s}^2+\|f_\ast\|_{\dot{H}^s}^2) - \|\mathcal{R}(v+f_\ast)\|_{L^2(\mathbb{S}^{n-1})}^2\\
& \leq \|\mathcal{R}\|^2\|v\|_{\dot{H}^s}^2-2\mathrm{Re}\langle\mathcal{R}v,\mathcal{R}f_\ast\rangle_{L^2}.
\end{align*}
For convenience we recall the definition $\mathcal{S}:=\mathcal{R}(-\Delta)^{-\frac{s}{2}}$ from Theorem \ref{cor-2}, treating $\mathcal{S}$ as an operator from $L^2(\R^n)$ to $L^2(\mathbb{S}^{n-1})$. By Theorem 1.1 of \cite{BMS}, $(-\Delta)^{\frac{s}{2}}f_\ast$ equals a constant multiple of $\mathcal{S}^\ast 1$, and by Theorem 2.1 of \cite{BMS} the operator $\mathcal{S}\mathcal{S}^\ast$ preserves the class of constant functions. Hence there are constants $c_1, c_2\neq 0$ such that
\[
\langle \mathcal{R}v,\mathcal{R}f_\ast\rangle_{L^2}=c_1\langle (-\Delta)^{\frac{s}{2}}v,\mathcal{S}^\ast 1\rangle_{L^2}=c_2\langle v,f_\ast\rangle_{\dot{H}^s}=0,
\]
as desired.
\endproof
Given the above result, a natural question arises concerning the inequalities \eqref{FS} and \eqref{HLS}; we recall that in this case $q=\frac{2(n-1)}{n-2s}$ and $s\in(\frac{1}{2},\frac{n}{2})$. For \eqref{FS} it is proved in \cite{CFW} that
\begin{equation}\label{CFW-eq2}
C_{\mathrm{FS}}(n,s)\|(-\Delta)^{\frac{s}{2}}G\|_2^2-\|G\|_{q}^2\lesssim\inf_{G_\ast\in M_{\mathrm{FS}}}\|(-\Delta)^{\frac{s}{2}}(G-G_\ast)\|_{2}^2,
\end{equation}
while for \eqref{HLS},
\begin{equation}\label{conj}
C_{\mathrm{HLS}}(n,s)\|f\|_{q'}^2-\left|\int_{\R^{2(n-1)}}\frac{f(x)\overline{f(y)}}{|x-y|^{n-2s}}\,\mathrm{d}x\mathrm{d}y\right|\lesssim \|f\|_{q'}^2\inf_{f_\ast\in M_{\mathrm{HLS}}}\left(\frac{\|f-f_\ast\|_{q'}}{\|f\|_{q'}}\right)^\sigma
\end{equation}
is proved in \cite{LZ} for $\sigma=1$ and is conjectured for $1<\sigma\leq 2$. An argument similar to the one used to prove Theorem \ref{t1-loc} shows that \eqref{conj} would follow from \eqref{CFW-eq2} if we knew that   
\begin{equation}\label{holder-near-rev2}
\|h_1-D_{r'}h_2\|_{r}^\sigma\gtrsim 1-\left(\int h_1 h_2\right)^2
\end{equation}
held for $r<2$ and unit vectors $h_1\in L^r$, $h_2\in L^{r'}$ such that $\int h_1h_2\geq 0$. When $\sigma=1$, \eqref{holder-near-rev2} is immediate from H\"older's inequality, and so we recover \eqref{conj} in this case. 

Although we do not know if \eqref{holder-near-rev2} is true for some $\sigma>1$ we have the following, which shows that it fails in the case of particular interest $\sigma=2$ and $r<2$.
\begin{prop}\label{HR2-fails}
If $1<r<\infty$, then a necessary condition for \eqref{holder-near-rev2} to hold is $\sigma\leq r$. In particular, if $r<2$ then \eqref{holder-near-rev2} is false for $\sigma=2$.
\end{prop}
\proof It is enough to disprove the weaker inequality
\begin{equation}\label{target}
\|h_1-D_{r'}h_2\|_{r}^\sigma\gtrsim 1-\int h_1 h_2.
\end{equation}
A result of Aldaz \cite{Aldaz} implies that
\begin{equation}\label{aldaz}
\big\||h_1|^{\frac{r}{2}}-|h_2|^{\frac{r'}{2}}\big\|_2^2\lesssim 1-\int |h_1h_2|.
\end{equation}
Combining \eqref{target} with \eqref{aldaz} we conclude that if $h_1$ and $h_2$ are non-negative then
\begin{equation}\label{target2}
\big\|h_1-h_2^{r'-1}\big\|_{r}^\sigma\gtrsim\big\|h_1^{\frac{r}{2}}-h_2^{\frac{r'}{2}}\big\|_2^2,
\end{equation}
and we claim that $\sigma\leq r$ is necessary for \eqref{target2} to hold. To prove this we modify an argument from \cite{CFL}: on the unit interval $[0,1]$, define
\[
h_1\equiv 1, \qquad\qquad h_2(x)=(1-\delta)^{-\frac{2}{r'}}\chi_{(0,(1-\delta)^2)}(x)
\] 
for $0<\delta\ll1$ fixed. Then clearly $\|h_1\|_{r}=\|h_2\|_{r'}=1$, and
\[
\big\|h_1^{\frac{r}{2}}-h_2^{\frac{r'}{2}}\big\|_2^2=(1-\delta)^2\left|1-\frac{1}{1-\delta}\right|^2+1-(1-\delta)^2=2\delta.
\]
Also,
\begin{align*}
\big\|h_1-h_2^{r'-1}\big\|_{r}^r & = (1-\delta)^2\left(\left|1-(1-\delta)^{-\frac{2}{r}}\right|^r+(1-\delta)^{-2}-1\right).
\end{align*}
By Taylor expansion, for this choice of $h_1$ and $h_2$, it follows that
\begin{equation*}
\big\|h_1-h_2^{r'-1}\big\|_{r}^r\lesssim\big\|h_1^{\frac{r}{2}}-h_2^{\frac{r'}{2}}\big\|_2^2,
\end{equation*}
and so by taking $\delta$ sufficiently small it follows that \eqref{target2} cannot hold whenever $\sigma>r$. \endproof

\begin{samepage}
\begin{remarks}\leavevmode\vspace{-1pt}
\begin{itemize}
\item The example we used in Proposition \ref{HR2-fails} may be easily generalised to $\R^n$, for example by letting $h_1$ be a positive radial Schwartz function of unit $L^r$ norm and taking $h_2(x):=(h_1(|x|))^{r-1}h\left(x'\right)$ (recall $x':=|x|^{-1}x$), where the function $h$ on $\mathbb{S}^{n-1}$ equals $(1-\delta)^{-\frac{2}{r'}}$ on a patch of measure $(1-\delta)^{2}$, and is zero otherwise.
\item The inequality \eqref{aldaz} may of course be viewed as a stable version of H\"older's inequality; such estimates have attracted attention in a range of contexts recently (see \cite{FO} and references therein). Although \eqref{aldaz} is weaker than \eqref{CFL-1} in general it still suffices for certain applications (as noted in \cite{CFL}), and in view of Proposition \ref{HR2-fails} it has the advantage that the reverse inequality is also true (see \cite{Aldaz}).
\end{itemize}
\end{remarks}
\end{samepage}
\begin{ack}
This work was supported by the JSPS Grant-in-Aid for Young Scientists A no.\ 16H05995 (Bez), by the JSPS KAKENHI 26287022 and 26610021 (Sugimoto), and by the JSPS Fellowship for Overseas Researchers FY 2015-2016 (Jeavons, Ozawa)
\end{ack}

\end{document}